\documentclass[10pt,a4paper]{article}
\usepackage{cite}
\usepackage{amscd}
\usepackage[latin1]{inputenc}
\usepackage{amsmath}
\usepackage{amsfonts}
\usepackage{amssymb}
\usepackage{color}
\usepackage{float}
\usepackage{amsthm}
\usepackage{breqn}

\usepackage{listings}
\newtheorem{theorem}{Theorem}

\newtheorem{proposition}[theorem]{Proposition}

\def\Qed{\hfill\raisebox{.6ex}{\framebox[2.5mm]{}}\\[.15in]}

\def\m{\mathbb}

\begin{document}
\date{}
\title{Cuspidal quintics and surfaces with\\ $p_g=0,$ $K^2=3$ and $5$-torsion}
\author{Carlos Rito}
\maketitle

\begin{abstract}

If $S$ is a quintic surface in $\mathbb P^3$ with singular set $15$ $3$-divisible ordinary cusps,
then there is a Galois triple cover $\phi:X\to S$ branched only at the cusps such that 
$p_g(X)=4,$ $q(X)=0,$ $K_X^2=15$ and $\phi$ is the canonical map of $X$.
We use computer algebra to search for such quintics having a free action of $\mathbb Z_5$,
so that $X/{\mathbb Z_5}$ is a smooth minimal surface of general type with $p_g=0$ and $K^2=3$.
We find two different quintics, one of which is the Van der Geer--Zagier quintic, the other is new.

We also construct a quintic threefold passing through the $15$ singular lines of the
Igusa quartic, with $15$ cuspidal lines there.
By taking tangent hyperplane sections, we compute quintic surfaces with singular set $17\mathsf A_2$, $16\mathsf A_2$, $15\mathsf A_2+\mathsf A_3$ and $15\mathsf A_2+\mathsf D_4$.


\end{abstract}

\section{Introduction}

In the context of Hilbert modular surfaces, Van der Geer and Zagier \cite{GeZa} have
constructed a smooth minimal surface of general type $X$ such that $K_X^2=10,$ $p_g(X)=4,$
$q(X)=0$ and the canonical map $\phi$ of $X$ is of degree $2$ onto a quintic surface $S$
with $20$ ordinary double points (nodes), ramified only over these points.
These two surfaces have a free action of $\m Z_5,$ thus by taking quotients one gets a
numerical Campedelli surface $X'$ ($K^2=2,$ $p_g=q=0$) which is a double cover of a
Godeaux surface $S'$ ($K^2=1,$ $p_g=q=0,$ $\pi_1=\m Z_5$), ramified over $4$ nodes.
Moreover, as noted by Catanese \cite[\S 5]{Ca}, $X$ is simply connected,
so also $\pi_1(X')=\m Z_5.$ This construction is a particular case of the Campedelli surfaces constructed in \cite{Ca}.

There is a similar construction involving a triple cover ramified over ordinary cusps
(singularities of type $\mathsf A_2$) instead of a double cover ramified over nodes.
The same paper \cite{GeZa} contains the construction of a quintic surface
$S$ in $\m P^3$ with $15$ cusps.
Later Barth \cite[unpublished]{Ba} has shown that this set of cusps is $3$-divisible, i.e.
there is a labelling $A_i,$ $A_i'$ for the $(-2)$-curves corresponding to the resolution
$A_i+A_i',$ $i=1,\ldots,15,$ of the cusps such that $$\sum_1^{15} (2A_i+A_i')\equiv 3L$$
for some divisor $L.$
Following Tan \cite[Thm 4.3.1]{Ta2} (see also \cite[\S 1.3]{Ta1}), this implies the existence
of a Galois triple cover $X\to S,$ ramified only over the cusps, such that $X$ is a smooth
minimal surface of general type with $K^2=15,$ $p_g=4,$ $q=0$ and this triple cover is the
canonical map of $X.$

The quintic $S$ has a free action of $\m Z_5,$ which extends to a free action on $X.$
So as noticed by Tan \cite[Theorem I]{Ta2},
taking quotients we get surfaces $X',$ $S'$ such that $K_{X'}^2=3,$ $p_g(X')=0,$
$K_{S'}^2=1,$ $p_g(S')=0$ and there is a triple cover $X'\rightarrow S'$ ramified only
over the $3$ cusps of $S'$.

The Van der Geer--Zagier quintic $S$ seems very special, it is invariant for the action
of the symmetric group in five elements.
In this paper we construct another quintic surface in $\m P^3$ with a $3$-divisible set
of $15$ cusps and with a free action of $\m Z_5$, and show that it is not isomorphic to $S$.

The observation which allowed us to have (computational) success with the construction is
the following. We considered some random quartic surfaces in $\m P^3$ with (at least) $15$
nodes at points $p_1,\ldots,p_{15}.$
For each case, the linear system of quintic surfaces with double points at $p_1,\ldots,p_{15}$
is of dimension $4,$ and we were able to find an element with $15$ cusps.
So our strategy is to first compute a $\m Z_5$-invariant quartic with $15$ nodes,
and then try to find the quintic.

Consider the $\m Z_5$ action $$(x:y:z:w)\mapsto (ex:e^2y:e^3z:e^4w),$$ where $e$ is a 5-th root of unity.
The quartic monomials invariant under this action are
\begin{displaymath}
x^3y,\ y^3w,\ xz^3,\ zw^3,\ x^2w^2,\ y^2z^2,\ xyzw.
\end{displaymath}
Thus a general invariant quartic passes with multiplicity $1$ through the four fixed points
$(1:0:0:0),\ldots,(0:0:0:1)$ of the action.
Semi-invariant quartics contain exactly three of these points and are singular at one of them.
For instance the quartic monomials corresponding to the eigenvalue 4 are
\begin{displaymath}
x^4,\ y^3z,\ yw^3,\ x^2zw,\ xy^2w,\ xyz^2,\ z^2w^2.
\end{displaymath}
A general quartic in this space passes through $(0:1:0:0), (0:0:0:1)$ and is singular at $(0:0:1:0)$.
So the semi-invariant quartic that we are going to compute has in fact $16$ nodes, while the invariant
quartic corresponding to the Van der Geer--Zagier quintic has exactly $15$ nodes. We use this to show that
our quintic surface is not isomorphic to the Van der Geer--Zagier quintic.

Now recall that the Igusa quartic threefold is singular at $15$ lines, so
by analogy with the quintic surface and its associated quartic, we might expect 
the existence of a quintic threefold in $\mathbb P^4$ with $15$ {\em cuspidal lines}, i.e. such that a general hyperplane section is a quintic surface with $15$ cusps. This is in fact the case, we compute an equation for such
threefold (in particular with $S_5$ symmetry). We think that this quintic is interesting on its own.

By taking hyperplanes tangent to this threefold, we found equations for quintic surfaces with
singular set $17\mathsf A_2$, $16\mathsf A_2$, $15\mathsf A_2+\mathsf A_3$ and $15\mathsf A_2+\mathsf D_4$.
To our knowledge only examples with at most $15$ cusps ($\mathsf A_2$ singularities) were known.

This suggests that quintic surfaces with more than $17$ cusps could exist.
The problem of constructing surfaces with many cusps has been addressed by
Barth and Rams (\cite{BaRa1}, \cite{BaRa2}).
The difficulty in finding new examples increases with the number of cusps.
Quintic surfaces with $n\geq 18$ cusps could provide interesting examples
of surfaces with non-birational canonical map (see e.g. \cite{Ta2}).

The paper is organized as follows.
In Section \ref{s2} we review the Van der Geer--Zagier quintic surface.
Then we explain the main steps of the construction of our quintic.
In section \ref{s4} we prove the $3$-divisibility of the $15$ cusps of our surface.
Then we compute the equation of the quintic threefold and finally
Section \ref{s6} explains the computation of the quintic surfaces with more than $15$ cusps.

We use the Computational Algebra System Magma \cite{BoCaPl} to perform the computations.
The corresponding input code lines are given in the Appendix (a file containing the output lines
is available on the author's webpage).

\bigskip
\noindent{\bf Notation}

We work over the complex numbers. All varieties are assumed to be projective algebraic.
A $(-n)$-curve on a surface is a curve isomorphic to $\m P^1$ with self-intersection $-n.$
Linear equivalence of divisors is denoted by $\equiv.$
The rest of the notation is standard in Algebraic Geometry.\\

\bigskip
\noindent{\bf Acknowledgements}

The author wishes to thank an anonymous Referee for his/her constructive comments.

The author is a member of the Center for Mathematics of the University of Porto.
This research was partially supported by FCT (Portugal) under the project PTDC/MAT-GEO/0675/2012 and by CMUP (UID/MAT/00144/2013), which is funded by FCT with national (MEC) and European structural funds through the programs FEDER, under the partnership agreement PT2020.

\section{Van der Geer--Zagier's quintic}\label{s2}

Van der Geer and Zagier \cite{GeZa} have shown the existence of a quintic surface $V$ in $\m P^3$ with $15$ ordinary cusps (singularities of type $\mathsf A_2$) and no other singularities.
It is given in $\m P^4(x,y,z,w,t)$ by $$12s_5-5s_2s_3=0,\ \ s_1=0,$$ where $s_i:=x^i+y^i+z^i+w^i+t^i.$
The $\m Z_5$ action $(x,y,z,w,t)\mapsto (y,z,w,t,x)$ acts freely on $V.$ The quotient of $V$ by this action is a Godeaux surface ($p_g=0,$ $K^2=1,$ $\pi_1=\m Z_5$) with $3$ cusps.

Later Barth \cite{Ba} showed that these cusps are $3$-divisible. He found a relation between the $(-2)$-curves contained in the resolution of the cusps and some divisors corresponding to $15$ lines contained in $V.$ We use an analogous argument below in Section \ref{s4}.

Here we note that exists exactly one quartic surface $Q$ with singular set $15$ ordinary
double points at the points $p_1,\ldots,p_{15}$ where $V$ has cusps.
Due to the symmetry, we expect it to be in the pencil generated by $s_4,$ $s_2^2,$
and in fact it is given in $\m P^4$ by $$4s_4-s_2^2=0,\ \ s_1=0.$$
In general it can be shown that $15$ nodes determine a quartic surface.

This surface contains $10$ plane conics such that each conic contains $6$ of the $p_i$'s. 
We show in the Appendix that $V$ and $Q$ meet exactly at these conics and confirm the
uniqueness of $Q.$

\section{The $\mathbb Z_5$-invariant surface}\label{s3}

\begin{theorem}
Let $q:=e^3+e^2,$ with $e$ a $5$-th root of unity.

\noindent
The $\m Z_5$-invariant quintic surface $S$ with equation
\begin{multline*}
3x^5 - 4y^5 + (60q + 120)xy^3z + (-90q - 150)x^2yz^2 + (-220q - 356)z^5 +\\
(30q + 30)x^2y^2w + (20q + 30)x^3zw + (390q + 630)yz^3w + (-210q - 350)y^2zw^2 +\\
(120q + 195)xz^2w^2 + (-120q - 180)xyw^3 + (20q + 32)w^5=0
\end{multline*}
has $15$ cusps and no other singularities.

\noindent
The quartic surface $Q$ with equation
\begin{multline*}
x^4 + (4q + 8)y^3z + (-12q - 20)xyz^2 + (4q + 4)xy^2w +\\
(4q + 6)x^2zw + (8q + 13)z^2w^2 + (-8q - 12)yw^3=0
\end{multline*}
has $16$ nodes, $15$ of which coincide with the cusps of $S$.

\noindent
Moreover, $S$ is not isomorphic to the Van der Geer--Zagier quintic.
\end{theorem}

\noindent{\bf Proof:}\\
The first two sentences are easily verified using computer algebra. This is done in the Appendix.
Suppose that $S$ is isomorphic to the Van der Geer--Zagier quintic $V.$
Then this isomorphism preserves the canonical systems of $S$ and $V,$
i.e. hyperplanes of $S$ are mapped to hyperplanes of $V,$ which implies that $S$ and $V$
are projectively equivalent. Hence the quartic $Q$ is mapped to
a $16$-nodal surface $Q'$ singular at the $15$ cusps of $V.$
But, as noted in Section \ref{s2}, the surface $Q'$ has exactly $15$ nodes.
\Qed

The surface $S$ is invariant for the action $$(x:y:z:w)\mapsto (x:ey:e^2z:e^3w)$$ and,
since its equation contains the monomials $x^5,y^5,z^5,w^5,$ the action is base point free.\\

Now we explain the steps taken to find the equation of $S$.
The corresponding Magma code is given in the Appendix.\\

We start by searching for the equation of a quartic surface $Q$ in $\m P^3$ with an action
of $\m Z_5$ and with $15$ ordinary double points.\\

Fix the $\m Z_5$ action as above.
The system of quartic polynomials which are invariant under this action is generated by
\begin{displaymath}
x^4,\ y^3z,\ yw^3,\ x^2zw,\ xy^2w,\ xyz^2,\ z^2w^2.
\end{displaymath}

We compute the monomials $s_1,s_2,s_3$ which generate the subsystem of the elements singular at the point $(1:1:1:1)$.\\

We want to find coefficients $U,V$ such that the polynomial $F:=s_1+Us_2+Vs_3$
defines a normal surface with $15$ nodes at the $\m Z_5$ orbits of points
$$p_1:=(1:1:1:1),\ p_2:=(x:y:z:w),\ p_3:=(a:b:c:d).$$
These are given by points in the scheme $S_C$ defined by
$F=\frac{\partial F}{\partial x}=\cdots =\frac{\partial F}{\partial w}=0$.\\

At this stage a generic point in $S_C$ corresponds to a non-normal quartic surface.
To overcome this, we impose the condition that the double point at $(1:1:1:1)$ is ordinary
(not all order $3$ minors of the Hessian matrix of $F$ vanish).\\

Still this does not give $15$ singular points. We need to add conditions to assure that
the points $p_1, p_2, p_3$ are in different $\m Z_5$ orbits:
$$(xd)^5-(aw)^5\ne 0,\ x^5-w^5\ne 0,\ a^5-d^5\ne 0.$$

Now the scheme $S_C$ is zero dimensional. Computing a point in $S_C,$ we get a $\m Z_5$-invariant quartic
surface $Q$ with $15$ double points $p_1,\ldots,p_{15}.$
As explained in the Introduction, there is an extra double point at $(0:0:1:0),$ a fixed point for the $\m Z_5$ action.\\

The linear system of $\m Z_5$-invariant polynomials of degree $5$ with double points at $p_1,\ldots,p_{15}$ has two generators $L_1, L_2.$ We want to compute $b$ such that the quintic surface with equation $F:=L_1+bL_2=0$ has cusps at the points $p_1,\ldots,p_{15}.$ This is done by imposing the vanishing of derivatives of $F$ and of all order $3$ minors of the Hessian matrix of $F$ at these points.\\

Finally we confirm, using Magma, that all singularities are of type $\mathsf A_2.$

\section{$3$-divisibility of the cusps}\label{s4}

Let $S$ be the $\m Z_5$-invariant quintic computed above and $G'$ be the quotient of $S$ by the $\m Z_5$ action.
Denote by $G$ the Godeaux surface obtained by resolving the three cusps of $G'.$

\begin{proposition}
The cusps of $G'$ are $3$-divisible, i.e. there exists a divisor $L$ such that $\sum_1^{3}(2A_i+A_i')\equiv 3L,$
for some labelling $A_i,$ $A_i'$ for the $(-2)$-curves in the resolution $A_i+A_i',$ $i=1,\ldots,3,$ of the cusps.
\end{proposition}
\noindent Notice that this implies the $3$-divisibility of the $15$ cusps of $S.$\\

\noindent{\bf Proof:}\\
Recall that a quartic Kummer surface has $16$ {\em tropes}, which are plane conics through $6$ nodes.
The quintic $S$ meets the Kummer surface $Q$ at the $10$ tropes of $Q$ not containing the 16th node $(0:0:1:0).$
These conics are divided into two $\m Z_5$ orbits $\overline T_1,$ $\overline T_2.$
One of the tropes of $Q$ ($\{y=0\}$) through $(0:0:1:0)$ is fixed by the $\m Z_5$ action and the remaining five
tropes through this point correspond to an orbit $\overline T_3$ of plane quintics in $S$ with five double points each.

Let $T_1,$ $T_2,$ $T_3$  be the strict transforms in $G$ of the quotients of $\overline T_1,$ $\overline T_2,$ $\overline T_3,$ respectively.
We compute the intersection matrix of $A_i, A_i',$ $i=1,2,3,$ and $T_1,$ $T_2,$ $T_3.$
This is achieved using the Magma function \verb|Blowup| to resolve the singularities of $S$ at a representative
of each of the three orbits of five cusps. The matrix is:

$$\bordermatrix{ & A_1 & A_1' & A_2 & A_2' & A_3 & A_3' & T_1 & T_2 & T_3 \cr
                A_1 & -2 & 1 & 0 & 0 & 0 & 0 & 1 & 1 & 2 \cr
                A_1' & 1 & -2 & 0 & 0 & 0 & 0 & 0 & 2 & 2 \cr
                A_2 & 0 & 0 & -2 & 1 & 0 & 0 & 0 & 2 & 1 \cr
                A_2' & 0 & 0 & 1 & -2 & 0 & 0 & 2 & 0 & 1 \cr
                A_3 & 0 & 0 & 0 & 0 & -2 & 1 & 1 & 1 & 2 \cr
                A_3' & 0 & 0 & 0 & 0 & 1 & -2 & 2 & 0 & 2 \cr                
                T_1 & 1 & 0 & 2 & 0 & 1 & 2 & -4 & 0 & 0 \cr
                T_2 & 1 & 2 & 0 & 2 & 1 & 0 & 0 & -4 & 0 \cr
                T_3 & 2 & 2 & 1 & 1 & 2 & 2 & 0 & 0 & -1 }$$\\
Its determinant is zero. Since $b_2(G)=9,$ these curves are dependent in ${\rm Num}(G)$ and this relation must be expressed in the nullspace of the matrix. As the basis for this nullspace is
$$\begin{pmatrix}
2 & 4 & 2 & -2 & -2 & -4 & -3 & 3 & 0
\end{pmatrix}$$
and $q(G)=0$ implies ${\rm NS}(G)={\rm Pic}(G)$ (Castelnuovo),
then there is a $5$-torsion element $t$ such that
$$2A_1+4A_1'+2A_2-2A_2'-2A_3-4A_3'-3T_1+3T_2\equiv t.$$
But $t\equiv t+5t\equiv 6t,$ therefore there exists $L$ such that
$$2A_1+A_1'+2A_2+A_2'+A_3+2A_3'\equiv 3L.$$

The Magma code for the above is in the Appendix.\Qed

\section{The quintic threefold}\label{s5}

It is well known that the maximum number of nodes that a cubic threefold can have is $10$, and there is exactly one such threefold,
the {\em Segre cubic}, which can be given in $\m P^4$ by the equation
$$x^3+y^3+z^3+w^3+t^3+h^3=0,\ \ h:=-x-y-z-w-t.$$
The dual of the Segre cubic is the so-called {\em Igusa quartic}. Its singular set is an union of $15$ lines.

Let the {\em monomial symmetric polynomial} $m_{(\alpha_1,\ldots,\alpha_5)}(x,y,z,w,t)$ be
defined as the sum of all monomials $x^{\alpha_i}\cdots t^{\alpha_j}$ over all distinct
permutations of $(\alpha_1,\ldots,\alpha_5)$.

\begin{proposition}
Let $X$ be the $S_5$-invariant quintic threefold with equation
$$7m_{(5,0,0,0,0)}-5m_{(4,1,0,0,0)}-2m_{(3,2,0,0,0)}+4m_{(3,1,1,0,0)}+$$
$$2m_{(2,2,1,0,0)}-4m_{(2,1,1,1,0)}+8m_{(1,1,1,1,1)}=0.$$
Then $X$
intersects the Igusa quartic along its singular lines, and has cusps along each of these lines.
\end{proposition}

This result is proved easily using computer algebra. Such an equation is computed as follows.\\

Fix $15$ points in the $15$ singular lines of the Igusa quartic and compute the linear system of $S_5$-invariant
quintics in $\mathbb P^4$ singular at these points.\\
We get a pencil of quintics which are singular at the $15$ lines.\\

Let $F, G$ be the generators of this pencil. Then compute $b$ such that the threefold $\{F+bG=0\}$
has cusps at these $15$ lines. This is done by imposing the vanishing of the order $3$ minors of the
Hessian matrix of $F+bG$ at the above $15$ points.\\

The corresponding Magma code is in the Appendix.

\section{Quintics with $17\mathsf A_2$, $16\mathsf A_2$, $15\mathsf A_2+\mathsf A_3$, $15\mathsf A_2+\mathsf D_4$}\label{s6}

A general hyperplane section of the quintic threefold computed above is a surface with $15$ ordinary cusps.
A tangent hyperplane section gives a surface $\{f=0\}$ with an extra double point,
which is ordinary only if one of the order $3$ minors of the Hessian matrix of $f$ is non-zero.
By imposing the vanishing of all order $3$ minors, we found quintic surfaces with singular set
$17\mathsf A_2$, $16\mathsf A_2$, $15\mathsf A_2+\mathsf A_3$ and $15\mathsf A_2+\mathsf D_4$.

The Appendix contains the computer code for the case $17\mathsf A_2$.
This surface is the hyperplane section $$3x-13y-13z+7w+7t=0$$ of the quintic computed on the previous section.

For the remaining cases the hyperplane is not defined over the rationals. The details can be found on the Author's webpage.

\appendix
\section*{Appendix: Magma code}\label{appendix}

The following code is implemented on the Computational Algebra System Magma, version V2.21-8.

\subsection*{Van der Geer--Zagier's quintic}

Here we verify the assertions made in Section \ref{s2} about the existence of a quartic surface $Q$ which is singular at the $15$ singular points
of the Van der Geer--Zagier quintic $V.$
\begin{verbatim}
P3<x,y,z,w>:=ProjectiveSpace(Rationals(),3);
t:=-(x+y+z+w);
s:=[x^i+y^i+z^i+w^i+t^i:i in [1..5]];
V:=Surface(P3,12*s[5]-5*s[2]*s[3]);
pts:=SingularPoints(V);#pts eq 15;
HasSingularPointsOverExtension(V) eq false;
L4:=LinearSystem(P3,4);
L:=LinearSystem(L4,[P3!x:x in pts],[2:x in pts]);
#Sections(L) eq 1;
Q:=Surface(P3,4*s[4]-s[2]^2);
Q eq Surface(P3,Sections(L)[1]);
#SingularPoints(Q) eq 15;
HasSingularPointsOverExtension(Q) eq false;
pc:=PrimeComponents(V meet Q);pc;
[#Points(pc[i] meet SingularSubscheme(V)) eq 6:i in [1..10]];
\end{verbatim}

\subsection*{The $\mathbb Z_5$-invariant surface}

\begin{verbatim}
K<e>:=CyclotomicField(5);
P<x,y,z,w>:=ProjectiveSpace(K,3);
\end{verbatim}
The following sequences of monomials $s_4$ and $s_5$ are invariant
under the action $(x:y:z:w)\mapsto (x:ey:e^2z:e^3w)$.
We compute the elements of $s_4$ which are singular at the point $(1:1:1:1)$.
\begin{verbatim}
s4:=[x^4,y^3*z,x*y*z^2,x*y^2*w,x^2*z*w,z^2*w^2,y*w^3];
s5:=[x^5,x^3*z*w,x^2*y^2*w,x^2*y*z^2,x*y^3*z,x*y*w^3,x*z^2*w^2,
     y^5,y^2*z*w^2,y*z^3*w,z^5,w^5];
LinearSystem(LinearSystem(P,s4),P![1,1,1,1],2);

A12<x,y,z,w,a,b,c,d,U,V,N1,N2>:=AffineSpace(K,12);
s:=[
    x^4 - 2*x^2*z*w + z^2*w^2,
    y^3*z - 2*x*y^2*w + x^2*z*w - z^2*w^2 + y*w^3,
    x*y*z^2 - x*y^2*w - z^2*w^2 + y*w^3
];
F:=s[1]+U*s[2]+V*s[3];
\end{verbatim}
We want to compute $U, V$ such that $\{F=0\}$ has double
points at (the orbits of) $(1:1:1:1)$, $(x:y:z:w)$ and $(a:b:c:d)$.
\begin{verbatim}
M:=Submatrix(HessianMatrix(Scheme(A12,F)),1,1,4,4);
m:=Minors(Evaluate(M,[1,1,1,1,a,b,c,d,U,V,N1,N2]),3)[1];
\end{verbatim}
If there is a number $N_1$ such that $1+mN_1=0$ at $(1:1:1:1)$,
then the minor $m\ne 0$ and the double point at $(1:1:1:1)$ is ordinary.
This reduces the probability of getting a non-normal surface.
\begin{verbatim}
Hx:=[F] cat [Derivative(F,i):i in [1..4]];
Ha:=[Evaluate(Hx[i],[a,b,c,d,a,b,c,d,U,V,N1,N2]):i in [1..#Hx]];
\end{verbatim}
To search for the $15$ double points, we define the scheme $S_C$ of the points
such that $H_x=H_a=0$ and $m\ne 0$.
\begin{verbatim}
SC:=Scheme(A12,Hx cat Ha cat [w-1,d-1,1+N1*m]);
\end{verbatim}
We impose extra conditions to ensure that the points
$(1:1:1:1)$, $(x:y:z:w)$ and $(a:b:c:d)$ are in different $\m Z_5$ orbits:
\begin{verbatim}
SC:=Scheme(SC,[1+N2*(x^5-1)*(x^5-a^5)*(a^5-1)]);
Dimension(SC) eq 0;
\end{verbatim}
We compute the points in $\m P^3$ corresponding to one of the points of $S_C$:
\begin{verbatim}
p:=Points(SC)[1];
p:=[P![1,1,1,1],P![p[1],p[2],p[3],p[4]],P![p[5],p[6],p[7],p[8]]];
\end{verbatim}
The unique $\m Z_5$-invariant quartic with nodes at (the
orbits of) the above three points:
\begin{verbatim}
L:=LinearSystem(LinearSystem(P,s4),p,[2,2,2]);
Q:=Scheme(P,Sections(L)[1]);
r:=SingularPoints(Q);
#r eq 16;
\end{verbatim}
The quartic $Q$ has $16$ nodes.
One of these points is $(0:0:1:0)$, a fixed point for the action of $\m Z_5$.
Now we compute the linear system of invariant quintics with nodes at
the $15$ points different from $(0:0:1:0)$:
\begin{verbatim}
L:=LinearSystem(LinearSystem(P,s5),p,[2,2,2]);
#Sections(L) eq 2;
\end{verbatim}
It remains to find the element of this pencil with $15$ cusps.
\begin{verbatim}
A5<X,Y,Z,W,b>:=AffineSpace(K,5);
h:=hom<CoordinateRing(P)->CoordinateRing(A5)|[X,Y,Z,W]>;
F:=h(Sections(L)[1])+b*h(Sections(L)[2]);
M:=HessianMatrix(Scheme(A5,F));RemoveColumn(~M,5);RemoveRow(~M,5);
H:=[Evaluate(Minors(M,3)[i],[1,1,1,1,b]):i in [1..#Minors(M,3)]];
\end{verbatim}
If all these minors vanish, the surface $\{F=0\}$ has a non-ordinary
double point at $(1:1:1:1)$.
\begin{verbatim}
SC:=Scheme(A5,H cat [X,Y,Z,W]);
Dimension(SC) eq 0;
\end{verbatim}
The points in $S_C$ give three possibilities for $b$. One of these corresponds to a
quintic surface with $15$ ordinary cusps ($\mathsf A_2$ singularities):
\begin{verbatim}
b:=Points(SC)[1][5];
F:=Sections(L)[1]+b*Sections(L)[2];
S:=Surface(P,F);
r:=SingularPoints(S);
\end{verbatim}
There are exactly $15$ singular points:
\begin{verbatim}
#r eq 15;
HasSingularPointsOverExtension(S) eq false;
\end{verbatim}
We confirm that these singularities are of type $\mathsf A_2$:
\begin{verbatim}
for x in r do IsSimpleSurfaceSingularity(S!x);end for;
\end{verbatim}

\subsection*{$3$-divisibility of the cusps}

The surfaces $S$ and $Q$ meet at the $10$ tropes of $Q$ not containing the point $(0:0:1:0)$:
\begin{verbatim}
pp:=PrimeComponents(Q meet S);
[#Points(pp[i] meet SingularSubscheme(S)):i in [1..10]];

P<x,y,z,w>:=P;
psi:=map<P->P|[x,y*e,z*e^2,w*e^3]>;
\end{verbatim}
These conics are divided into two $\m Z_5$ orbits $\overline T_1$, $\overline T_2$:
\begin{verbatim}
T1:=&join[pp[i]:i in [1,2,5,6,9]];
T1 eq psi(T1);
T2:=&join[pp[i]:i in [3,4,7,8,10]];
T2 eq psi(T2);
[Multiplicity(T1,r[i]):i in [4,5,6]] eq [1,2,3];
[Multiplicity(T2,r[i]):i in [4,5,6]] eq [3,2,1];
\end{verbatim}
We define $\overline T_3$ as the $\m Z_5$ orbit of plane quintics corresponding to $5$ of the $6$
tropes of $Q$ through the fixed point $(0:0:1:0)$:
\begin{verbatim}
L1:=LinearSystem(LinearSystem(P,1),[P!r[i]:i in [1,3,11,12,13]]);
t:=[Scheme(P,Sections(L1)[1])];
#Points(t[1] meet SingularSubscheme(S)) eq 5;
for i in [1..4] do t:=t cat [psi(t[#t])];end for;
T3:=(&join t) meet S;
T3 eq psi(T3);
[Multiplicity(T3,r[i]):i in [4,5,6]] eq [4,2,4];
\end{verbatim}
The strict transform of $T_3$ does not intersect the strict transform of $T_1+T_2$:
\begin{verbatim}
Points(T3 meet (T1 join T2)) eq r;
\end{verbatim}
The resolution of each cusp is an union of two $(-2)$-curves $A_i, A_i'$.
We compute the intersection number of $T_1$, $T_2$ and $T_3$ with these curves.
First we blowup the quintic at a cusp:
\begin{verbatim}
T:=[T1,T2,T3];
for j in [4,5,6] do
 X,mp:=Blowup(S,P!r[j]);
\end{verbatim}
We define the exceptional divisor $E$ and the $(-2)$-curves $A_i, A_i'$:
\begin{verbatim}
 E:=(P!r[j]) @@ mp;
 A:=PrimeComponents(E meet Complement(X,E));
\end{verbatim}
and compute the intersection number of the strict transform of $T_i$ with the $(-2)$-curves:
\begin{verbatim}
 for i in [1,2,3] do
  t:=Complement(Scheme(S,DefiningEquations(T[i])) @@ mp,E);
  #Points(t meet A[1]),#Points(t meet A[2]);
 end for;
end for;
\end{verbatim}
This shows that the intersection matrix of the curves, in the Godeaux
surface, corresponding to the $\mathsf A_2$ configurations of the three cusps and
the orbits $T_1$, $T_2$ and $T_3$ is:
\begin{verbatim}
M:=SymmetricMatrix([-2,1,-2,0,0,-2,0,0,1,-2,0,0,0,0,-2,0,0,
0,0,1,-2,1,0,0,2,1,2,-4,1,2,2,0,1,0,0,-4,2,2,1,1,2,2,0,0,-1
]);M;
Determinant(M) eq 0;
Nullspace(M);
\end{verbatim}

\subsection*{The quintic threefold}

We define the Segre cubic $S_3$, the Igusa quartic $I_4$
and compute the singular set of $I_4$:
\begin{verbatim}
P4<x,y,z,w,t>:=ProjectiveSpace(Rationals(),4);
S3:=Scheme(P4,x^3+y^3+z^3+w^3+t^3+(-x-y-z-w-t)^3);
rho:=map<P4->P4|Basis(JacobianIdeal(DefiningEquation(S3)))>;
I4:=rho(S3);
SI4:=SingularSubscheme(I4);
\end{verbatim}
The $S_5$ symmetric quintics are generated by:
\begin{verbatim}
e:=[ElementarySymmetricPolynomial(CoordinateRing(P4),i):i in [1..5]];
s5:=[e[1]^5, e[1]^3*e[2], e[1]^2*e[3], e[1]*e[4],
e[1]*e[2]^2, e[2]*e[3], e[5]];
\end{verbatim}
We fix $15$ singular points of $I_4$ and compute the symmetric quintics singular at these points:
\begin{verbatim}
r:=Points(Scheme(SI4,x+2*y+3*z+4*w+5*t));                                          
L5:=LinearSystem(P4,s5);
L:=LinearSystem(L5,[P4!x:x in r],[2:x in r]);
\end{verbatim}
A generic element of the pencil $L$ is singular exactly at
the $15$ singular lines of the Igusa:
\begin{verbatim}
X1:=SingularSubscheme(Scheme(P4,Sections(L)[1]));
X2:=SingularSubscheme(Scheme(P4,Sections(L)[2]));
SI4 eq ReducedSubscheme(X1);
SI4 eq ReducedSubscheme(X2);
\end{verbatim}
We compute the element of $L$ which contains $15$ cuspidal lines:
\begin{verbatim}
A<X,Y,Z,W,T,b>:=AffineSpace(Rationals(),6);
h:=hom<CoordinateRing(P4)->CoordinateRing(A)|[X,Y,Z,W,T]>;
F:=h(Sections(L)[1])+b*h(Sections(L)[2]);
\end{verbatim}
We compute $b$ such that the quintic $\{F = 0\}$ has cusps
at the $15$ points above: 
\begin{verbatim}
H:=HessianMatrix(Scheme(A,F));RemoveColumn(~H,6);RemoveRow(~H,6);
G:=[F] cat [Derivative(F,i):i in [1..5]] cat Minors(H,3);
\end{verbatim}
The vanishing of these minors implies non-ordinary double points. 
\begin{verbatim}
G:=[Evaluate(G[j],Coordinates(r[i]) cat [b]):
j in [1..#G],i in [1..#r]];
S:=Scheme(A,G cat [X,Y,Z,W,T]);
Dimension(S) eq 0;
PointsOverSplittingField(S);
\end{verbatim}
This gives $b = -5/7$.\\\\
The quintic threefold with $15$ cuspidal lines:
\begin{verbatim}
F:=Sections(L)[1]-5/7*Sections(L)[2];
Q:=Scheme(P4,F);
SI4 eq ReducedSubscheme(SingularSubscheme(Q));
Degree(SingularSubscheme(Q)) eq 30;
\end{verbatim}
The coefficients of $F$ in ${\rm Sections}(L_5)$:
\begin{verbatim}
CoefficientMap(L5)(F);
\end{verbatim}
We verify that a random hyperplane section is a quintic surface with $15$ cusps:
\begin{verbatim}
a:=[Random(1,100):i in [1..5]];
S:=Surface(P4,[F,a[1]*x+a[2]*y+a[3]*z+a[4]*w+a[5]*t]);
#SingularPoints(S) eq 15;
r:=SingularPoints(S);
for x in r do IsSimpleSurfaceSingularity(S!x);end for;
\end{verbatim}

\subsection*{The quintic surface with $17$ cusps}

Now we compute an hyperplane section of the above quintic threefold
which is a quintic surface with $17$ cusps ($\mathsf A_2$ singularities):
\begin{verbatim}
A<X,Y,Z,W,T,a,b,c,d>:=AffineSpace(Rationals(),9);
h:=hom<CoordinateRing(P4)->CoordinateRing(A)|[X,Y,Z,W,T]>;
G:=a*X+b*Y+c*Z+d*W+T;
f:=Evaluate(h(F),T,-(a*X+b*Y+c*Z+d*W));
\end{verbatim}
So $\{f=0\}$ is an hyperplane section of the quintic threefold.
\begin{verbatim}
H:=HessianMatrix(Scheme(A,f));
H:=Submatrix(H,1,1,5,5);
min1:=Minors(H,3);
\end{verbatim}
If all these minors vanish, then the double point is non-ordinary.
\begin{verbatim}
J:=[JacobianSequence(h(F))[i]:i in [1..5]];
min2:=Minors(Matrix([J,[a,b,c,d,1]]),2);
\end{verbatim}
If all these minors vanish, then the hyperplane is tangent to the quintic.
\begin{verbatim}
H:=[T-1,h(F),G] cat min1 cat min2;
SC:=Scheme(A,H cat [X,Y-2]);
pp:=PrimeComponents(SC);
[Dimension(x):x in pp];
\end{verbatim}
Since $S_C$ has zero dimensional components, we can compute some points.
These data allow us to compute the quintic surface with $17$ $\mathsf A_2$ points:
\begin{verbatim}
p:=Points(pp[4]);
a:=p[1][6];b:=p[1][7];c:=p[1][8];d:=p[1][9];
Q:=Surface(P4,[F,a*x+b*y+c*z+d*w+t]);
r:=Points(SingularSubscheme(Q));
#r eq 17;
for i in [1..#r] do IsSimpleSurfaceSingularity(Q!r[i]);end for;
HasSingularPointsOverExtension(Q);
\end{verbatim}

\bibliography{ReferencesRito}

\begin{thebibliography}{BoCaPl}
\expandafter\ifx\csname urlstyle\endcsname\relax
  \expandafter\ifx\csname doi\endcsname\relax
  \def\doi#1{doi:\discretionary{}{}{}#1}\fi \else
  \expandafter\ifx\csname doi\endcsname\relax
  \def\doi{doi:\discretionary{}{}{}\begingroup \urlstyle{rm}\Url}\fi \fi

\bibitem[Ba]{Ba}
W.~Barth, {\em A quintic surface with 15 three-divisible Cusps\/} (2000),
  {P}reprint, Erlangen.

\bibitem[BaRa1]{BaRa1}
W.~Barth and S.~Rams, {\em Equations of low-degree projective surfaces with
  three-divisible sets of cusps\/}, Math. Z., {\bf 249} (2005), no.~2,
  283--295.

\bibitem[BaRa2]{BaRa2}
W.~Barth and S.~Rams, {\em Cusps and codes\/}, Math. Nachr., {\bf 280} (2007),
  no. 1-2, 50--59.

\bibitem[BoCaPl]{BoCaPl}
W.~Bosma, J.~Cannon and C.~Playoust, {\em The {M}agma algebra system. {I}.
  {T}he user language\/}, J. Symbolic Comput., {\bf 24} (1997), no. 3-4,
  235--265, ISSN 0747-7171, computational algebra and number theory (London,
  1993).

\bibitem[Ca]{Ca}
F.~Catanese, {\em {Babbage's conjecture, contact of surfaces, symmetric
  determinantal varieties and applications}\/}, Invent. Math., {\bf 63} (1981),
  433--465.

\bibitem[Ta1]{Ta1}
S.-L. Tan, {\em Surfaces whose canonical maps are of odd degrees\/}, Math.
  Ann., {\bf 292} (1992), no.~1, 13--29.

\bibitem[Ta2]{Ta2}
S.-L. Tan, {\em Cusps on some algebraic surfaces and plane curves\/}, {\em
  Complex Analysis, Complex Geometry and Related Topics - Namba\/}, vol.~60,
  2003, 106--121.

\bibitem[GeZa]{GeZa}
G.~van~der Geer and D.~Zagier, {\em {The Hilbert modular group for the field
  $\mathbb Q(\sqrt{13})$}\/}, Invent. Math., {\bf 42} (1977), 93--133.

\end{thebibliography}

\

\

\noindent Carlos Rito\\
\\{\it Permanent address:}
\\ Universidade de Tr\'as-os-Montes e Alto Douro, UTAD
\\ Quinta de Prados
\\ 5000-801 Vila Real, Portugal
\\ www.utad.pt
\\ crito@utad.pt\\
\\{\it Current address:}
\\ Departamento de Matem\' atica
\\ Faculdade de Ci\^encias da Universidade do Porto
\\ Rua do Campo Alegre 687
\\ 4169-007 Porto, Portugal
\\ www.fc.up.pt
\\ crito@fc.up.pt

\end{document}